# Integral Triangles with one angle twice another,

# and with the bisector (of the double angle) also of integral length.

*Konstantine Zelator*
*P.O. Box 4280*
*Pittsburgh, PA 15203*
*U.S.A*
*e-mail address:* *konstantine_zelator@yahoo.com*

## 1. Introduction

The purpose of this work is two-folded. First, to parametrically describe all integral triangles with one angle being twice another. The second goal is to also parametrically describe a subset of the aforementioned integral triangles. Namely, those integral triangles that not only have the property that one of their angles is twice another; but with the additional feature that the bisector of the double angle is also of integral length. (Recall that an integral triangle is one with all three side lengths being integers).

Author M.N Deshpande, in an article published in the *Mathematical Gazette* in 2002 (see Reference **[1]**), showed that some integral triangles with angle *B* being twice angle *A* ; are generated by the formulas $a = n^2$ , $b = mn$ , $c = m^2 - n^2$; where *m* and *n* are positive integers such that $n < m < 2n$; and *a,b,c* are the side lengths opposite from angles *A, B, C* respectively.

As we show in **Result 1** (found in Section 5 of this paper), the entire family of integral triangles with angle *B* being twice angle *A* can be described by the following parametric formulas: $a = lk^2$, $b = lkm$, $c = l(m^2 - k^2)$, where $l, k, m$ are positive integers such that $m$ and $k$ are relative prime; and with either $k < m < k\sqrt{2}$ (i.e. $k^2 < m^2 < 2k^2$) or alternatively $k\sqrt{2} < m < 2k$ (i.e. $2k^2 < m^2 < 4k^2$) . These then are precise. (i.e. necessary and sufficient) conditions for describing the entire set of integral triangles having one angle being twice another.

In Section 6, we identify the subfamily (of the integral triangles with one angle being twice another) in which every triangle has angle bisector (of the double angle) also of integral length. This is done in **Result 2** wherein we give a precise parametric description of those triangles. All triangles, whether integral or not, with angle *B* being twice angle *A*; must satisfy the condition, $b^2 = a(a+c)$ where $a, b, c$ are the three side lengths. We prove that this condition is necessary and sufficient for a triangle to have angle *B* twice angle *A*. This is done in Proposition 1, in Section 4 of this article. This same condition can also be found in a paper by Wynne William Wilson, published in the *Mathematical Gazette* in 1976 (see Reference **[2]**) .

In **Lemma 2** (Section 3) we show that if $a,b,c$ are three positive real numbers such that $a(a+c)=b^2$. Then, a triangle having $a,b,c$ as its side lengths; can be formed if and only if either $0<c\leq a$; or alternatively $0<a<c<3a$ . Such a triangle will have according to Proposition 1, angle *B* being twice angle *A* .

**Lemmas 4, 5, 6** in Section 3, are well known in number theory; and can easily be found in number theory books. For instance, see References **[3]** and **[4]**.

In Section 2, we present five pictures of triangles with one angle being twice another. These five pictures describe the five general settings, in which a triangle has angle *B* twice angle *A*. We finish this paper by offering some closing remarks in Section 7.

**Section 2: Triangles with one angle being twice another: Five Figures**

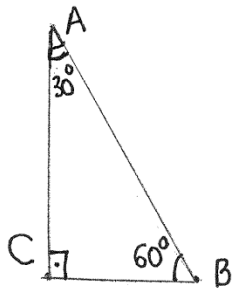

$30°-60°-90°$ triangle

$\angle B = 60°$ , $\angle A = 30°$

**Figure 1**

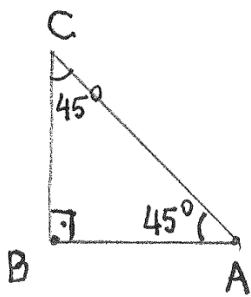

$45°-45°-90°$ triangle

$\angle B = 90°$ , $\angle A = 45°$

**Figure 2**

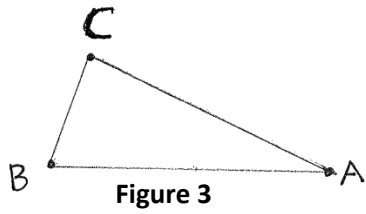

Acute triangle

$\angle B = 2\angle A$

**Figure 3**

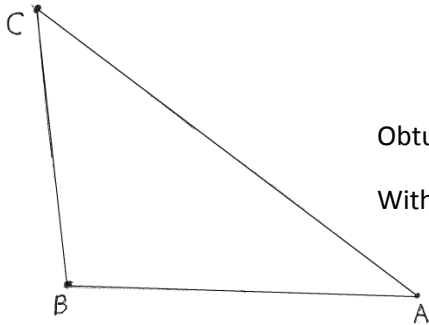

Obtuse triangle

With the double angle $B$ being the obtuse angle $\angle B = 2\angle A$

**Figure 4**

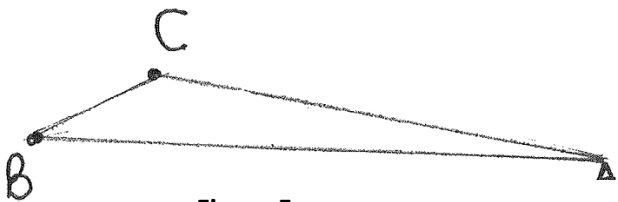

Obtuse triangle with one of the two

acute angles being twice the other $\angle B = 2\angle A$

(Obtuse angle at $C$)

**Figure 5**

**Section 3: Six Lemmas**

**Lemma 1:** *Let $ABC$ be a triangle, and* c, b, a *the side lengths of the sides $\overline{AB}$, $\overline{AC}$, and $\overline{CB}$ respectively. Then $c = a\cos B + b\cos A$.*

**Proof:** From the Law of Cosines we have,

$$a^2 = b^2 + c^2 - 2bc\cos A$$
$$b^2 = c^2 + a^2 - 2ca\cos B$$ **(1).**

Adding the two equations in **(1)** member wise yields,

$$a^2 + b^2 = 2c^2 + a^2 + b^2 - 2c(b\cos A + a\cos B);$$
$$2c2 = 2c(b\cos A + a\cos B); \text{ or equivalently (since } c > 0)$$
$$c = b\cos A + a\cos B \;\square$$

**Lemma 2:** *Let a,b,c be positive real numbers such that* $b^2 = a(a+c)$. *Then a, b, c are the side lengths of a triangle if and only if; either* $0 < c \leq a$ *; or alternatively* $0 < a < c < 3a$.

**Proof**: Three positive reals $a, b, c$ are the side lengths of a triangle if and only if, the three triangle inequalities are satisfied: $a < b+c$, $b < a+c$, $c < a+b$. So let $a, b, c$ be three positive real numbers such that $\qquad b^2 = a(a+c) \qquad$ **(2).**

First observe that two of the three triangles inequalities above follow from condition **(2)** alone. These two inequalities are $a < b+c$ and $b < a+c$. Indeed, since $a$ and $c$ are positive, $a+c > a$; and so $(a+c)^2 > a(a+c)$; and by **(2)**, we infer that $(a+c)^2 > b^2$; which implies $|a+c| > |b|$; and so (since $a > 0$, $b > 0$, $c > 0$) we obtain $a+c > b$.

Now the other inequality $a < b+c$.

We have $(b+c)^2 = b^2 + 2bc + c^2 = a(a+c) + 2bc + c^2 = a^2 + ac + 2bc + c^2 > a^2$, since $a, b, c$ are positive by **(2)**. Therefore, $(b+c)^2 > a^2$ from which it follows that $b+c > a$.

Now the third inequality, $c < a+b.$

We will show that under **(2)**; if $0 < c \leq a$ or $0 < a < c < 3a$. Then $c < a+b$. And conversely, if $c < a+b$ and **(2)** hold, then either $0 < c \leq a$ or $0 < a < c < 3a$.

Suppose that $0 < c \leq a$ or $0 < a < c < 3a$. If $0 < c \leq a$; then since $b > 0$, we have $0 < c < c+b$. **(3)**.

But also (from $0 < c \leq a$ ), $0 < b < c+b \leq a+b$. **(4)**.

From **(3)** and **(4)** it follows that $0 < c < a+b$.

Now suppose that $0 < a < c < 3a$ **(5)**.

From **(5)**, it follows that

$$c(c-3a) < 0 \Leftrightarrow c^2 - 3ac < 0 \Leftrightarrow c^2 - 2ac + a^2 < a^2 + ac = a(a+c); \quad \textbf{(6)}.$$

By **(6)** and **(2)** it follows that, $(c-a)^2 < b^2$; $|c-a| < |b|$.

And since $c > a$ (by **(5)**) and $b > 0$. The last inequality gives $c - a < b$; $c < a+b$.

Next, we show the converse: namely that if $0 < c < a+b$; then either $0 < c \leq a$ or $0 < a < c < 3a$.

From $0 < c < a+b$ it follows that, $c - a < b$ **(7)**.

If $c \leq a$, then $0 < c \leq a$ (since c, a, b are positive), then we are done. If on the other hand $a < c$; then $0 < c - a$. And so by **(7)**, we have $0 < c - a < b$ **(8)**.

The inequalities **(8)** imply $(c-a)^2 < b^2 \Leftrightarrow c^2 - 2ac + a^2 < b^2$ (and by **(2)**) we further get,

$c^2 - 2ac + a^2 < a(a+c)$; $c^2 - 2ac + a^2 < a^2 + ac \Leftrightarrow c(c-3a) < 0$ which in view of **(8)**,

$a > 0$ and $c > 0$; implies $0 < a < c < 3a$.

The proof is complete. □

**Lemma 3:** *Let $ABC$ be a triangle and $\overline{BD}$ the angle bisector of the angle $\angle B$. Also, let a,b,c be the lengths of the sides $\overline{BC}$, $\overline{CA}$, and $\overline{AB}$; respectively. Also, let $|AD|$ and $|DC|$ denote the lengths of the line segments $\overline{AD}$ and $\overline{DC}$, respectively. Then $|AD| = \dfrac{cb}{a+b}$ and $|DC| = \dfrac{ab}{a+c}$.*

**Proof:**

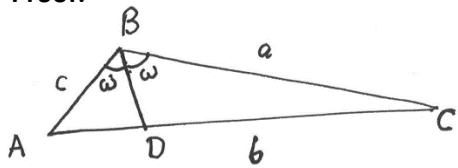

**Figure 6**

By the well-known theorem in Euclidean Geometry, we have

$$\frac{|AD|}{|DC|} = \frac{c}{a} \qquad (9).$$

From **(9)** we obtain,

$$\frac{|AD| + |DC|}{|DC|} = \frac{c+a}{a};$$

$$\frac{|AC|}{|DC|} = \frac{c+a}{a};$$

And so,

$$\frac{b}{|\overline{DC}|} = \frac{c+a}{a}; \text{ which implies}$$

$$|\overline{DC}| = \frac{ab}{c+a}.$$

Also from **(9)**,

$$\frac{|\overline{DC}|}{|\overline{AD}|} = \frac{a}{c};$$

$$\frac{|\overline{DC}| + |\overline{AD}|}{|\overline{AD}|} = \frac{a+c}{c};$$

$$\frac{b}{|\overline{AD}|} = \frac{a+c}{c}; \text{ which gives}$$

$$|\overline{AD}| = \frac{bc}{a+c}.$$

The proof is complete. ☐

For Lemma 4 below, also known as Euclid's Lemma, see references **[3]** and **[4].**

**Lemma 4:** *Let* a,b,c *be positive integers, with* a *and* b *being relatively prime. In addition, suppose that* a *is a divisor of the product* bc. *Then the integer must be a divisor of the integer* c.

For the two lemmas below, **Lemmas 5** and **6**; the interested reader may refer to **[4].**

**Lemma 5:** *Suppose that* a *and* b *are positive integers and* n *a natural number. Furthermore, assume that the integer* $a^n$ *is a divisor of the integer* $b^n$. *Then, the integer* a *is a divisor of the integer* b.

**Lemma 6:** *Let* n *be a natural number; and* a,b,c *positive integers such that* $ab=c^n$, *and with* a *and* b *being relatively prime. Then there exists positive integers* $a_1$ *and* $b_1$, *such that* $a=a^n_1$, $b=b^n_1$, $c=a_1 b_1$, *and with* $a_1$ *and* $b_1$ *being relatively prime.*

**4. A Key Condition**

The following proposition gives a necessary and sufficient condition for a triangle to contain an angle which is twice another angle.

**Proposition 1:** *Let ABC be a triangle with* a,b,c *being the side lengths of the sides* $\overline{BC}$, $\overline{CA}$, *and* $\overline{AB}$ *respectively. Then the angle at B is twice the angle at A if, and only if,* $b^2 = a(a+c)$.

**Proof:** First suppose that $B = 2(angleA)$ or simply $b = 2A$ **(10).**

We will show that **(10)** implies $b^2 = a(a+c)$. From the Law of Sines we have

$$2R = \frac{a}{\sin A} = \frac{b}{\sin B} = \frac{c}{\sin C} \quad \textbf{(11)},$$

Where *R* is the radius of the circumscribed circle. Using the double angle identity $\sin 2A = 2\sin A \cos A$ **(10)**, and **(11)** we get

$$\frac{a}{\sin A} = \frac{b}{2\sin A \cos A}; \text{ which implies}$$

$$\cos A = \frac{b}{2a} \quad \textbf{(12)}.$$

From **(10)** and **Lemma 1** we further have, $c = b\cos A + a\cos 2A$      **(13).**

And by **(12), (13),** and the double-angle identity $\cos 2A = 2\cos^2 A - 1$ we further obtain

$$c = b\left(\frac{b}{2a}\right) + a\left[2\left(\frac{b}{2a}\right)^2 - 1\right];$$

$$c = \frac{b^2}{2a} + \frac{b^2}{2a} - a;$$

$$2ac = 2b^2 - 2a^2;$$

$$b^2 = a(a+c).$$

Now the converse statement, assume that $b^2 - a(a+c)$ **(14).**

By **Lemma 1** and **(14)** we have $b^2 = a(a + b\cos A + a\cos B)$ **(15).**

Combining **(11)** with **(15)** yields,

$$4R^2 \sin^2 B = 4R^2 \sin A[\sin B\cos A + \sin A\cos B + \sin A] \,;$$

$\sin^2 B = \sin A[\sin B\cos A + \sin A\cos B + \sin A]$ **(16).**

From the summation identity $\sin(A+B) = \sin B\cos A + \sin A\cos B$ and **(16)** we get,

$$\sin^2 B = \sin A[\sin(A+B) + \sin A] \,;$$

$$\sin^2 B - \sin^2 A = \sin A \sin(A+B) \,;$$

$(\sin B - \sin A)(\sin B + \sin A) = \sin A \sin(A+B)$      **(17).**

Combining **(17)** with the summation identities

$$\sin B - \sin A = 2\sin(\frac{B-A}{2})\cos(\frac{B+A}{2}) \text{ and } \sin B + \sin A = 2\sin(\frac{A+B}{2})\cos(\frac{A-B}{2}) \text{ ; further gives}$$

$$4\sin(\frac{B-A}{2})\cos(\frac{B+A}{2})\sin(\frac{A+B}{2})\cos(\frac{A-B}{2}) = \sin A \sin(A+B) \textbf{ (18)}.$$

Using the fact that $A+B+C=180°$ and consequently $\sin C = \sin(A+B)$ and also that

$$\sin(\frac{C}{2}) = \cos(\frac{B+A}{2}) \text{ and } \sin(\frac{A+B}{2}) = \cos(\frac{C}{2}) \text{ ; \textbf{(18)} gives}$$

$$[2\sin(\frac{B-A}{2})\cos(\frac{A-B}{2})][2\cos(\frac{C}{2})\sin(\frac{C}{2})] = \sin A \sin C \textbf{ (19)}.$$

From the double-angle identity for the sine and the fact that

$$\cos(\frac{A-B}{2}) = \cos(\frac{B-A}{2}) \text{ ; \textbf{(19)} further yields}$$

$\sin(B-A)\sin C = \sin A \sin C$ ; and since $\sin C \neq 0$ (in fact $\sin C > 0$, since C is a triangle angle) we get $\sin(B-A) = \sin A$ **(20)**.

Since $A$ and $B$ are triangle angles; they fall strictly between 0 and 180 degrees. So, since $\sin A > 0$, **(20)** implies that $\sin(B-A) > 0$. And so the difference $B-A$ must also between 0 and 180 degrees. Hence by **(20)**, either $A = B-A$ or the absolute value of the difference $(B-A)-A$ is equal to 180 degrees; which means either $B-2A = 180°$; or $A-(B-A) = 180°$. We cannot have $B-2A = 180°$; $B = 2A + 180°$ since B cannot exceed 180 degrees. Now if $A-(B-A) = 180°$ ;

$2A = 180° + B; A = 90° + \dfrac{B}{2}$ , which means that $A$ is strictly between 90 and 180 degrees; while both triangle angles $B$ and $C$ must be acute angles.

But this cannot occur since then $B - A = B - (90° + \dfrac{B}{2}) = \dfrac{B}{2} - 90° < 0°$ and also $-180° < B - A$ (since $A$ and $B$ are triangle angles). Thus $-180° < B - A < 0°$ , which implies that $\sin(B-A) < 0$ , contradicting **(20)** (since $\sin A > 0$ ).

Going back to **(20)** we see that **(20)** implies $B - A = A$ ; or equivalently $B = 2A$ . The proof is complete $\square$

The following proposition is an immediate consequence of **Lemma 2** and **Proposition 1**.

**Proposition 2:** *A triangle* $ABC$ *can be formed with the angle at* $B$ *being twice the angle at* $A$ *if, and only if,* $b^2 = a(a+c)$ *and either* $0 < c \leq a$ *or alternatively* $0 < a < c < 3a$ ; *where a,b,c are the side lengths of the sides* $\overline{BC}, \overline{CA}$ *, and* $\overline{AB}$ *respectively.*

**5. Integral triangle with one angle being twice another.**

Consider an integral triangle ABC with the angle B being twice the angle at A. Then according to Proposition 2, the positive integers a,b,c must satisfy the precise conditions $b^2 = a(a+c)$ and either $0 < c \leq a$ (and so $1 \leq c \leq a$ , since c is a positive integer) or alternatively $0 < a < 3a$ (and so $1 \leq a < c < 3a$ ) **(21).**

Let $\ell$ be the greatest common divisor of *a* and *c*. Then, $a = \ell a_1 , c = \ell c_1$ where $\ell, a_1, c_1$ are positive integers and with $a_1$ and $c_1$ being relatively prime **(22).**

From **(22)** and **(21)** we obtain, $b^2 = \ell^2 \cdot a_1 \cdot (a_1 + c_1)$ **(23).**

According to **(23)**, $\ell^2$ is a divisor of $b^2$. Thus, by **Lemma 5**, it follows that $\ell$ must be a divisor of b. And so, $b = \ell b_1$, for some positive integer $b_1$ **(24).**

Combining **(24)** with **(23)** gives, $b_1^2 = a_1(a_1 + c_1)$ **(25).**

Since $a_1$ and $c_1$ are relatively prime; so must be $a_1$ and $a_1 + c_1$. Thus by **(25)** and **Lemma 6** it follows that $a_1 = k^2$, $a_1 + c_1 = m^2$, $b_1 = km$.

For positive integers *k,m* with $k < m$, and *k,m* being relatively prime. Moreover, from **(24)** and **(22)** we also get $b = \ell km$, $a = \ell k^2$, $c = \ell(m^2 - k^2)$

Applying the inequality conditions of **(21)** we arrive at the integer conditions $k^2 < m^2 \leq 2k^2$ or alternatively $2k^2 < m^2 < 4k^2$. Note that the equal sign cannot hold, since m, k are integers and the $\sqrt{2}$ is irrational.

We can now state **Result 1**.

**Result 1:** *The entire family of integral triangles* ABC, *with angle at* B *being twice the angle at* A; *can be parametrically described (in terms of three integer parameters) as follows:*

*Side lengths* $a = \ell k^2, b = \ell km, c = \ell(m^2 - k^2)$ ; *where* $\ell, m, k$ *are positive integers such that k and m are relatively prime and with either,* $k^2 < m^2 < 2k^2$ *or alternatively,* $2k^2 < m^2 < 4k^2$.

## 6. Integral triangles with one angle being twice another and with the bisector (of the double angle) also of integral length.

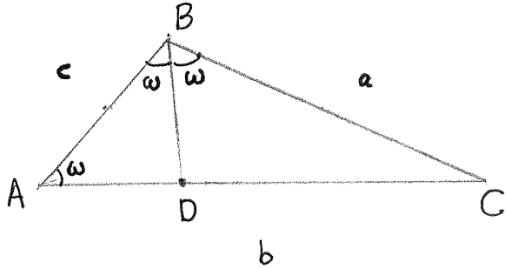

**Figure 7**

In Figure 6, a triangle *ABC* is illustrated with angle at $B = 2\omega = 2$ times the angle at *A*. Let $\overline{BD}$ be the angle bisector of the angle at *B*. We have $<CAB = \omega =< ABD =< DBC$ .

The triangle *ABD* is isosceles with the apex at *D*. We set $|AD| = r = |BD|$ **(26).**

And suppose that the triangle *ABC* is an integral triangle with *r* also an integer. Then by **Result 1, Lemma 3**, and **(26)**; we have

$$b = \ell km, a = \ell k^2, c = \ell(m^2 - k^2), |\overline{DC}| = \frac{\ell k^3}{m}, r = \frac{\ell(m^2 - k^2)k}{m}$$ ; where *k* and *m* are relatively prime

positive integers, $\ell$ of a positive integer; and such that either $k^2 < m^2 \leq 2k^2$ or alternatively

$2k^2 < m^2 < 4k^2$ **(27).**

From $r = \frac{\ell(m^2 - k^2)k}{m}$ in **(27)**; we see that since *r* is an integer *m* must be a divisor of $\ell$. Indeed, since

*m* and *k* are relatively prime; it follows that *m* and $m^2 - k^2$ must also be relatively prime. Thus the integer *m* must be relatively prime to the product $\ell \cdot \left[(m^2 - k^2)k\right]$ And so, (because *r* is an integer); it follows by Euclids's Lemma (**Lemma 4**) that m must be a divisor of $\ell$.

Hence, $\ell = m \cdot d$, For some positive integer *d* **(28).**

By combining **(28)** with **(27)** we obtain formulas that parametrically describe the integral triangles with the angle bisector of the double-angle being also integral.

**Result 2:** *The entire family of integral triangles* ABC *with angle* $\angle B = 2(\angle A)$ *and with the angle bisector* $\overline{BD}$ *being of integral length* $r = |\overline{BD}|$ *as well; can be parametrically described as follows:*

$b = dkm^2$, $a = dmk^2$, c=dm(m² − k²), r=$|\overline{BD}| = |\overline{AD}| = dk(m^2 - k^2)$, $|\overline{DC}| = dk^3$; *where* d,k,m *are positive integers with* k *and* m *being relatively prime; and with either* $k^2 < m^2 \leq 2k^2$ *or alternatively,* $2k^2 < m^2 < 4k^2$.

**7. Closing Remarks**

    **Remark 1:** Every member or triangle of the family of triangles described in **Result 2** has integral values for five of its lengths: *a, c, b,* $|\overline{BD}| = r = |\overline{AD}|$, and $|\overline{DC}|$; are all integers. Thus, triangle $BDC$ in **Figure 7** is an integral triangle. But it also has one angle being twice the other. Indeed, $\angle BDC = 2\omega$ and $\angle DBC = \omega$.

    **Remark 2:** Looking at the parametric formulas in **Result 2**, we see that the *primitive triangles* in the family described in **Result 2** are precisely those that satisfy $m = d = 1$. Indeed, a *primitive triangle* is one with the greatest common of the three sidelengths being 1. In the case of **Result 2,** we have $\gcd(a,b,c) = dm$, in virtue of the fact that $\gcd(k,m) = 1$. Thus $\gcd(a,b,c) = 1 \Leftrightarrow dm = 1$; $d = 1 = m$.